# On Some Algebra and the Corresponding Analysis in the Pseudo-Riemannian Space with Signature (1, -1, -1, -1).

D. M. Volokitin

### Introduction.

In this paper the certain 4-dimensional algebra in 4-dimensional pseudo-Riemannian space with signature (1, -1, -1, -1) is constructed. On the basis of this algebra the elements of the analysis, i.e. the theory of 4-dimensional functions of the 4-dimensional variable are built up. In the process of designing the analysis two additional assumptions about the properties of functions are made. Obtained under different assumptions about the properties of functions, the Cauchy-Riemann equations are solved in a flat, spherically and cylindrically symmetric cases.

### Algebra and analysis.

In order to regulate the interpretation I allow myself to repeat the main results of the article [7].

Let's consider the four-dimensional pseudo-Riemannian space with a metric $g_{ij}(x^0, x^1, x^2, x^3)$ that satisfies the following conditions:

$$g_{ij} = g_{ji}$$

$$g_{00} > 0$$

$$\begin{vmatrix} g_{00} & g_{01} \\ g_{10} & g_{11} \end{vmatrix} < 0$$

$$\begin{vmatrix} g_{00} & g_{01} & g_{02} \\ g_{10} & g_{11} & g_{12} \\ g_{20} & g_{21} & g_{22} \end{vmatrix} > 0$$

$$\begin{vmatrix} g_{00} & g_{01} & g_{02} & g_{03} \\ g_{10} & g_{11} & g_{12} & g_{13} \\ g_{20} & g_{21} & g_{22} & g_{23} \\ g_{30} & g_{31} & g_{32} & g_{33} \end{vmatrix} < 0$$

As it is known [3, p.301], the metric tensor in any frame of reference feasible with real bodies must satisfy these conditions.

Let's fix the point $(x_0^0, x_0^1, x_0^2, x_0^3)$ and metric $g_{ij}(x_0^0, x_0^1, x_0^2, x_0^3) = g_{ij}$ in it.

In the future, wherever it would be a flat space, in order to avoid cluttering the formulas we omit the subscript in expressions such as "$x_0^i$", pointing to the location.

Let's consider pseudo-Euclidean space tangent to original curved space in this point. We introduce all over the tangent space the metric $g_{ij}$.

For the above conditions imposed on the metric, there is a unique triangular transformation of the basis $\{\vec{e}_i\}$, resulting a quadratic form $\sum_{i,j=0}^{3} g_{ij} x^i x^j$ to have a normal form [2, p.195-199]

$$(y^0)^2 - (y^1)^2 - (y^2)^2 - (y^3)^2.$$

Further in all cases, where it will not cause misunderstandings, we omit the limits of summation. We'll find the matrix

$$A = (a_i^j) = \begin{pmatrix} a_0^0 & 0 & 0 & 0 \\ a_1^0 & a_1^1 & 0 & 0 \\ a_2^0 & a_2^1 & a_2^2 & 0 \\ a_3^0 & a_3^1 & a_3^2 & a_3^3 \end{pmatrix}$$

of the triangular transformation mentioned above $\vec{e}_i = \sum_j a_i^j \vec{k}_j$.

Instead of using the Jacobi method [1, p.195-199] of finding the matrix $A$ we use the fact that during the transformation of the basis $\{\vec{e}_i\}$ to the basis $\{\vec{k}_i\}$ with the help of the matrix $A$ the matrix

$$\eta = \begin{pmatrix} 1 & 0 & 0 & 0 \\ 0 & -1 & 0 & 0 \\ 0 & 0 & -1 & 0 \\ 0 & 0 & 0 & -1 \end{pmatrix}$$

is transformed according to the law $A \eta A^T$, where $A^T$ is the transposed matrix $A$ [1, str.189]. That is, we solve the system of equations

$$A \eta A^T = g \quad \text{or} \quad \sum_{m,n} a_i^m \eta_{mn} (a^T)_n^j = \sum_{m,n} \eta_{mn} a_i^m a_j^n = g_{ij}.$$

This system of non-linear equations in the variables $a_j^i$ due to the specific form of the matrix $A$ is easily solved by consistent exception of groups of variables $(a_0^0, a_1^0, a_2^0, a_3^0), (a_1^1, a_2^1, a_3^1), (a_2^2, a_3^2), a_3^3$.

Since the explicit expressions for $a_j^i$ are far too cumbersome, we shall not represent them here.

Further let be $\vec{k}_i = \sum_j b_i^j \vec{e}_j$, i.e. $B = (b_j^i) = A^{-1}$.

We introduce in the tangent space the rules of algebra multiplication of basis vectors $\vec{k}_i$

$$\vec{k}_0 \vec{k}_0 = \vec{k}_1 \vec{k}_1 = \vec{k}_2 \vec{k}_2 = \vec{k}_3 \vec{k}_3 = a\vec{R}$$
$$\vec{k}_1 \vec{k}_0 = \vec{k}_0 \vec{k}_1 = -\alpha^1 a\vec{R}$$
$$\vec{k}_2 \vec{k}_0 = \vec{k}_0 \vec{k}_2 = -\alpha^2 a\vec{R}$$
$$\vec{k}_3 \vec{k}_0 = \vec{k}_0 \vec{k}_3 = -\alpha^3 a\vec{R} \quad , \tag{1}$$
$$\vec{k}_1 \vec{k}_2 = -\vec{k}_2 \vec{k}_1 = \alpha^3 c\vec{R}$$
$$\vec{k}_2 \vec{k}_3 = -\vec{k}_3 \vec{k}_2 = \alpha^1 c\vec{R}$$
$$\vec{k}_3 \vec{k}_1 = -\vec{k}_1 \vec{k}_3 = \alpha^2 c\vec{R}$$

where $a, c, \alpha^0, \alpha^1, \alpha^2, \alpha^3$ are the real constants, and the conditions

$$\alpha^0 = 1, \; (\alpha^1)^2 + (\alpha^2)^2 + (\alpha^3)^2 = 1, \tag{2}$$
$$\vec{R} = \alpha^0 \vec{k}_0 + \alpha^1 \vec{k}_1 + \alpha^2 \vec{k}_2 + \alpha^3 \vec{k}_3 \tag{3}$$

are carried out.

Thus only two of three variables $\alpha^i$ are independent.

It is easy to verify the conditions $\vec{z}\vec{R} = \vec{R}\vec{y} = 0$ ($\vec{y}$ and $\vec{z}$ are random vectors) and $(\vec{e}_i \vec{e}_j)\vec{e}_k = \vec{e}_i(\vec{e}_j \vec{e}_k) = 0$ for any $i, j, k = 0,...,3$. Naturally, because of the last relationships, this algebra is associative.

In addition, the product of any element $y^0 \vec{k}_0 + y^1 \vec{k}_1 + y^2 \vec{k}_2 + y^3 \vec{k}_3$ on its conjugate $y^0 \vec{k}_0 - y^1 \vec{k}_1 - y^2 \vec{k}_2 - y^3 \vec{k}_3$ is as given

$$(y^0 \vec{k}_0 + y^1 \vec{k}_1 + y^2 \vec{k}_2 + y^3 \vec{k}_3)(y^0 \vec{k}_0 - y^1 \vec{k}_1 - y^2 \vec{k}_2 - y^3 \vec{k}_3) = ((y^0)^2 - (y^1)^2 - (y^2)^2 - (y^2)^2)\vec{R} \qquad (4)$$

with the factor before $\vec{R}$ having a shape typical to the pseudo-norm of Minkowski space.
Actually, these two were the main considerations during the designing of the given algebra.
It should also be noted that a linear transformation of the basis that preserves the form (4) of the product of any element on its conjugate, is a Lorentz transformation.
It is understood that the product of any three elements in this algebra is zero.

We call a right quotient of dividing the element $\vec{r}$ by the element $\vec{\rho}$ the element $\vec{\eta}$: $\{\dfrac{\vec{r}}{\vec{\rho}}\}_r = \vec{\eta}$, if $\vec{r} = \vec{\eta}\vec{\rho}$.

But $\vec{\eta}\vec{\rho} \sim \vec{R}$ (if $\vec{\rho} \neq \vec{R}$), that is why $\vec{r} \sim \vec{R}$.

Obviously by the arbitrary elements, except $\vec{R}$, can be divided from the right elements proportional to $\vec{R}$ only.

Performing the multiplication in the last equation, we get:
$\vec{r} = r^0 \vec{e}_0 + r^1 \vec{e}_1 + r^2 \vec{e}_2 + r^3 \vec{e}_3 =$
$= \vec{\eta}\vec{\rho} = (\eta^0 \vec{e}_0 + \eta^1 \vec{e}_1 + \eta^2 \vec{e}_2 + \eta^3 \vec{e}_3)(\rho^0 \vec{e}_0 + \rho^1 \vec{e}_1 + \rho^2 \vec{e}_2 + \rho^3 \vec{e}_3) =$
$= [(\eta^0 \rho^0 + \eta^1 \rho^1 + \eta^2 \rho^2 + \eta^3 \rho^3)a + (-\eta^1 \rho^0 - \eta^0 \rho^1)\alpha^1 a +$
$+ (-\eta^2 \rho^0 - \eta^0 \rho^2)\alpha^2 a + (-\eta^3 \rho^0 - \eta^0 \rho^3)\alpha^3 a + (-\eta^2 \rho^1 + \eta^1 \rho^2)\alpha^3 c +$
$+ (\eta^3 \rho^1 - \eta^1 \rho^3)\alpha^2 c + (-\eta^3 \rho^2 + \eta^2 \rho^3)\alpha^1 a](\vec{e}_0 + \alpha^1 \vec{e}_1 + \alpha^2 \vec{e}_2 + \alpha^3 \vec{e}_3)$

If we denote the expression in square brackets as $\tilde{A}$, the following relations must hold
$\tilde{A} = r^0$
$\tilde{A}\alpha^1 = r^1$
$\tilde{A}\alpha^2 = r^2$
$\tilde{A}\alpha^3 = r^3$

Right quotient can not be determined uniquely, since from one equation $\tilde{A} = r^0$ (the other three are proportional to the first equation) one can not determine four unknown $\eta^i$.

If $\vec{\rho} \sim \vec{R}$, then $\vec{r} = \vec{\eta}\vec{\rho} \sim \vec{\eta}\vec{R} = 0$, i.e. zero is the only element that can be divided from the right by elements proportional to $\vec{R}$.

All of the above is quite similar if we define the left quotient of the elements: $\{\dfrac{\vec{r}}{\vec{\rho}}\}_l = \vec{\eta}$, if $\vec{r} = \vec{\rho}\vec{\eta}$.

Thus, all the elements of algebra are zero divisors.

Now as we know the matrices $(a_i^j)$ and $(b_i^j)$ we can find the rules of multiplication of basis vectors $\vec{e}_i$:

$$\vec{e}_i \vec{e}_j = h_{ij} \vec{R} , \qquad (5)$$

where $\vec{R} = \sum_n \beta^n \vec{e}_n$, $\qquad (6)$

$$\beta^i = \sum_j \alpha^j b^i_j \qquad (7)$$

(here $(\alpha^i)$ is the constant isotropic vector $(1, \alpha^1, \alpha^2, \alpha^3)$ ),

$$h_{ij} = \sum_{m,n} a^m_i a^n_j \kappa_{mn}, \qquad (8)$$

$$(\kappa_{ij}) = \begin{pmatrix} a & -\alpha^1 a & -\alpha^2 a & -\alpha^3 a \\ -\alpha^1 a & a & \alpha^3 c & -\alpha^2 c \\ -\alpha^2 a & -\alpha^3 c & a & \alpha^1 c \\ -\alpha^3 a & \alpha^2 c & -\alpha^1 c & a \end{pmatrix}. \qquad (9)$$

It should be noted that the determinant of the last matrix is 0, the rank is 3.

Vector $(\beta^i)$ (or $\vec{R}$), of course, will also be isotropic.

Thus, in any point of the original pseudo-Riemannian space, you can enter the local algebra
$$\vec{e}_i(x^0, x^1, x^2, x^3) \vec{e}_j(x^0, x^1, x^2, x^3) = h_{ij}(x^0, x^1, x^2, x^3; \alpha^1, \alpha^2, \alpha^3, a, c) \vec{R}$$
where $\vec{R} = \sum_i \beta^i(x^0, x^1, x^2, x^3; \alpha^1, \alpha^2, \alpha^3, a, c) \vec{e}_i(x^0, x^1, x^2, x^3)$. But the parameters $a, c, \alpha^1, \alpha^2, \alpha^3$ are no longer constants, but are the functions of the coordinates of the point.

We introduce in the whole pseudo-Riemannian space the local geometry
$$d\vec{r}(x^0, x^1, x^2, x^3) = dx^0 \vec{e}_0(x^0, x^1, x^2, x^3) + dx^1 \vec{e}_1(x^0, x^1, x^2, x^3) + \\ + dx^2 \vec{e}_2(x^0, x^1, x^2, x^3) + dx^3 \vec{e}_3(x^0, x^1, x^2, x^3)$$

Let's introduce a variety of functions
$$\vec{u}(\vec{r}) = \{u^0(x^0, x^1, x^2, x^3), u^1(x^0, x^1, x^2, x^3), u^2(x^0, x^1, x^2, x^3), u^3(x^0, x^1, x^2, x^3)\}$$
on the manifold $\{x^i\}$, and on this variety the same local geometry and algebra are introduced:
$$D\vec{u}(x^0, x^1, x^2, x^3) = Du^0(x^0, x^1, x^2, x^3) \vec{e}_0(x^0, x^1, x^2, x^3) + \\ + Du^1(x^0, x^1, x^2, x^3) \vec{e}_1(x^0, x^1, x^2, x^3) + \\ + Du^2(x^0, x^1, x^2, x^3) \vec{e}_2(x^0, x^1, x^2, x^3) + \\ + Du^3(x^0, x^1, x^2, x^3) \vec{e}_3(x^0, x^1, x^2, x^3)$$
where
$$Du^i(x^0, x^1, x^2, x^3) = \sum_j [u^i_{,j}(x^0, x^1, x^2, x^3) + \sum_k \Gamma^i_{kj}(x^0, x^1, x^2, x^3) u^k(x^0, x^1, x^2, x^3)] dx^j = \\ = \sum_j u^i_{;j}(x^0, x^1, x^2, x^3) dx^j$$

Here $\Gamma^i_{kj}$ are Christoffel symbols and $u^i_{,j}$ and $u^i_{;j}$ are the usual and covariant derivatives of the vector components, respectively.

We call the left-hand derivative of the function $\vec{u}(\vec{r})$ in the point $\vec{r}$ not depending on the method of approach to zero the limit of the difference ratio $[\frac{D\vec{u}(\vec{r})}{d\vec{r}}]_l$:

$$\vec{u}'_l(\vec{r}) = \lim_{d\vec{r} \to 0} [\frac{D\vec{u}(\vec{r})}{d\vec{r}}]_l, \text{ if the division is performed from the left.}$$

Omitting the arguments of $Du^i(x^0, x^1, x^2, x^3)$, we get:

$$\vec{u}'_l(\vec{r}) = \lim_{d\vec{r} \to 0}[\frac{D\vec{u}(\vec{r})}{d\vec{r}}]_l =$$

$$= \lim_{\substack{dx^0 \to 0 \\ dx^1 \to 0 \\ dx^2 \to 0 \\ dx^3 \to 0}} \frac{Du^0\vec{e}_0 + Du^1\vec{e}_1 + Du^2\vec{e}_2 + Du^3\vec{e}_3}{dx^0\vec{e}_0 + dx^1\vec{e}_1 + dx^2\vec{e}_2 + dx^3\vec{e}_3} =$$

$$= \lim_{\substack{dx^0 \to 0 \\ dx^1 \to 0 \\ dx^2 \to 0 \\ dx^3 \to 0}} \{[(u^0_{;0}dx^0 + u^0_{;1}dx^1 + u^0_{;2}dx^2 + u^0_{;3}dx^3)\vec{e}_0 +$$

$$+ (u^1_{;0}dx^0 + u^1_{;1}dx^1 + u^1_{;2}dx^2 + u^1_{;3}dx^3)\vec{e}_1 + \quad . \tag{10}$$

$$+ (u^2_{;0}dx^0 + u^2_{;1}dx^1 + u^2_{;2}dx^2 + u^2_{;3}dx^3)\vec{e}_2 +$$

$$+ (u^3_{;0}dx^0 + u^3_{;1}dx^1 + u^3_{;2}dx^2 + u^3_{;3}dx^3)\vec{e}_3]/$$

$$/[dx^0\vec{e}_0 + dx^1\vec{e}_1 + dx^2\vec{e}_2 + dx^3\vec{e}_3]\}$$

We set $d\vec{r} = dx^i\vec{e}_i$. Then

$$\vec{u}'_l(\vec{r}) = \frac{u^0_{;i}\vec{e}_0 + u^1_{;i}\vec{e}_1 + u^2_{;i}\vec{e}_2 + u^3_{;i}\vec{e}_3}{\vec{e}_i} = v^0\vec{e}_0 + v^1\vec{e}_1 + v^2\vec{e}_2 + v^3\vec{e}_3 ,$$

where $v^i$ are the components of the left-hand derivative.
Hence we obtain

$$\sum_j u^j_{;i}\vec{e}_j = \sum_j v^j\vec{e}_i\vec{e}_j = \sum_j v^j h_{ij}\vec{R} = (\sum_j v^j h_{ij})\sum_k \beta^k\vec{e}_k = s_i\sum_j \beta^j\vec{e}_j ,$$

where $s_i = \sum_k v^j h_{ij}$. \tag{11.1}

and $u^j_{;i} = s_i\beta^j$, \tag{12.1}

We write the last relations in more detailed form.

$u^0_{;i} = s_i\beta^0$

$u^1_{;i} = s_i\beta^1$

$u^2_{;i} = s_i\beta^2$

$u^3_{;i} = s_i\beta^3$ .

Next, we carry out all possible divisions between left and right sides of the last equations. We obtain

$$\frac{u^0_{;i}}{u^1_{;i}} = \frac{\beta^0}{\beta^1} \qquad \frac{u^0_{;i}}{u^2_{;i}} = \frac{\beta^0}{\beta^2} \qquad \frac{u^0_{;i}}{u^3_{;i}} = \frac{\beta^0}{\beta^3}$$

$$\frac{u^1_{;i}}{u^0_{;i}} = \frac{\beta^1}{\beta^0} \qquad \frac{u^1_{;i}}{u^2_{;i}} = \frac{\beta^1}{\beta^2} \qquad \frac{u^1_{;i}}{u^3_{;i}} = \frac{\beta^1}{\beta^3}$$

$$\frac{u^2_{;i}}{u^0_{;i}} = \frac{\beta^2}{\beta^0} \qquad \frac{u^2_{;i}}{u^1_{;i}} = \frac{\beta^2}{\beta^1} \qquad \frac{u^2_{;i}}{u^3_{;i}} = \frac{\beta^2}{\beta^3}$$

$$\frac{u^3_{;i}}{u^0_{;i}} = \frac{\beta^3}{\beta^0} \qquad \frac{u^3_{;i}}{u^1_{;i}} = \frac{\beta^3}{\beta^1} \qquad \frac{u^3_{;i}}{u^2_{;i}} = \frac{\beta^3}{\beta^2} \quad .$$

It follows that

$$u^0_{;i}\beta^1 = u^1_{;i}\beta^0 \qquad u^0_{;i}\beta^2 = u^2_{;i}\beta^0 \qquad u^0_{;i}\beta^3 = u^3_{;i}\beta^0$$
$$u^1_{;i}\beta^0 = u^0_{;i}\beta^1 \qquad u^1_{;i}\beta^2 = u^2_{;i}\beta^1 \qquad u^1_{;i}\beta^3 = u^3_{;i}\beta^1$$
$$u^2_{;i}\beta^0 = u^0_{;i}\beta^2 \qquad u^2_{;i}\beta^1 = u^1_{;i}\beta^2 \qquad u^2_{;i}\beta^3 = u^3_{;i}\beta^2$$
$$u^3_{;i}\beta^0 = u^0_{;i}\beta^3 \qquad u^3_{;i}\beta^1 = u^1_{;i}\beta^3 \qquad u^3_{;i}\beta^2 = u^2_{;i}\beta^3 \ .$$

Adding here trivial expressions like $u^0_{;i}\beta^0 = u^0_{;i}\beta^0$, we get invariant relatively to coordinate transformations analogue of the Cauchy-Riemann equations:
$$u^j_{;i}\beta^k = u^k_{;i}\beta^j, \quad 0 \le i, j, k \le 3 \ .$$

Clearly, not all of these conditions are independent. We choose as independent, such as these:
$$u^j_{;i}\beta^0 = u^0_{;i}\beta^j, \quad 0 \le i \le 3, \quad 1 \le j \le 3 \ .$$

Now if we make
$$d\vec{r} = dx^0\vec{e}_0 + dx^1\vec{e}_1 + dx^2\vec{e}_2 + dx^3\vec{e}_3 = d\tau\vec{R} =$$
$$= d\tau\beta^0\vec{e}_0 + d\tau\beta^1\vec{e}_1 + d\tau\beta^2\vec{e}_2 + d\tau\beta^3\vec{e}_3 \ ,$$

(that is, take the isotropic direction of tending $d\vec{r}$ to zero), then
$$dx^i = \beta^i d\tau \ .$$

And further
$$\vec{u}'_l(\vec{r}) = \lim_{d\tau \to 0}\{[(u^0_{;0}d\tau\beta^0 + u^0_{;1}d\tau\beta^1 + u^0_{;2}d\tau\beta^2 + u^0_{;3}d\tau\beta^3)\vec{e}_0 +$$
$$+ (u^1_{;0}d\tau\beta^0 + u^1_{;1}d\tau\beta^1 + u^1_{;2}d\tau\beta^2 + u^1_{;3}d\tau\beta^3)\vec{e}_1 +$$
$$+ (u^2_{;0}d\tau\beta^0 + u^2_{;1}d\tau\beta^1 + u^2_{;2}d\tau\beta^2 + u^2_{;3}d\tau\beta^3)\vec{e}_2 + \ .$$
$$+ (u^3_{;0}d\tau\beta^0 + u^3_{;1}d\tau\beta^1 + u^3_{;2}d\tau\beta^2 + u^3_{;3}d\tau\beta^3)\vec{e}_3]/$$
$$/[d\tau(\beta^0\vec{e}_0 + \beta^1\vec{e}_1 + \beta^2\vec{e}_2 + \beta^3\vec{e}_3)]\}$$

From here we obtain $u^i_{;0}\beta^0 + u^i_{;1}\beta^1 + u^i_{;2}\beta^2 + u^i_{;3}\beta^3 = 0$, because zero is the only element that can be divided by element proportional to $\vec{R}$.

And, finally, the analogue of the Cauchy-Riemann equations:
$$\sum_j u^i_{;j}\beta^j = 0 \ , \quad 0 \le i \le 3$$
$$\beta^0 u^i_{;j} = \beta^i u^0_{;j}, \quad 1 \le i \le 3, \quad 0 \le j \le 3$$
(13)

However, not all of these equations are invariant under coordinate transformations. In the invariant record, they look like this:
$$\sum_j u^i_{;j}\beta^j = 0 \ , \quad 0 \le i \le 3$$
$$u^j_{;i}\beta^k = u^k_{;i}\beta^j, \quad 0 \le i, j, k \le 3$$
(13a)

For greater clarity, we write out in detail the conditions (13).
$$u^0_{;0}\beta^0 + u^0_{;1}\beta^1 + u^0_{;2}\beta^2 + u^0_{;3}\beta^3 = 0 \qquad (13.1)$$
$$u^1_{;0}\beta^0 + u^1_{;1}\beta^1 + u^1_{;2}\beta^2 + u^1_{;3}\beta^3 = 0 \qquad (13.2)$$
$$u^2_{;0}\beta^0 + u^2_{;1}\beta^1 + u^2_{;2}\beta^2 + u^2_{;3}\beta^3 = 0 \qquad (13.3)$$
$$u^3_{;0}\beta^0 + u^3_{;1}\beta^1 + u^3_{;2}\beta^2 + u^3_{;3}\beta^3 = 0 \qquad (13.4)$$

$$u^1_{;0} = \frac{\beta^1}{\beta^0} u^0_{;0} \qquad (13.5)$$

$$u^2_{;0} = \frac{\beta^2}{\beta^0} u^0_{;0} \tag{13.6}$$

$$u^3_{;0} = \frac{\beta^3}{\beta^0} u^0_{;0} \tag{13.7}$$

$$u^1_{;1} = \frac{\beta^1}{\beta^0} u^0_{;1} \tag{13.8}$$

$$u^2_{;1} = \frac{\beta^2}{\beta^0} u^0_{;1} \tag{13.9}$$

$$u^3_{;1} = \frac{\beta^3}{\beta^0} u^0_{;1} \tag{13.10}$$

$$u^1_{;2} = \frac{\beta^1}{\beta^0} u^0_{;2} \tag{13.11}$$

$$u^2_{;2} = \frac{\beta^2}{\beta^0} u^0_{;2} \tag{13.12}$$

$$u^3_{;2} = \frac{\beta^3}{\beta^0} u^0_{;2} \tag{13.13}$$

$$u^1_{;3} = \frac{\beta^1}{\beta^0} u^0_{;3} \tag{13.14}$$

$$u^2_{;3} = \frac{\beta^2}{\beta^0} u^0_{;3} \tag{13.15}$$

$$u^3_{;3} = \frac{\beta^3}{\beta^0} u^0_{;3} \tag{13.16}$$

If we take, for example, equations (13.5), (13.8), (13.11) and (13.14), multiply their left and right sides by $\beta^0$, $\beta^1$, $\beta^2$ and $\beta^3$, accordingly, then add them up term by term, we get:

$$u^1_{;0}\beta^0 + u^1_{;1}\beta^1 + u^1_{;2}\beta^2 + u^1_{;3}\beta^3 = \frac{\beta^1}{\beta^0}(u^0_{;0}\beta^0 + u^0_{;1}\beta^1 + u^0_{;2}\beta^2 + u^0_{;3}\beta^3).$$

Now, notice that the left side of this equation is the left-hand side of equation (13.2), and the right is multiplied by $\frac{\beta^1}{\beta^0}$ the left-hand side of equation (13.1). Thus, equation (13.2) is a consequence of (13.1). The same can be said about the equations (13.3) and (13.4).

So, only any one of equations (13.1-4) and equations (13.5-16) are independent.

Having thirteen equations for sixteen $u^i_{;j}$ unknown one can express the remaining $u^i_{;j}$ by any other three of them, for example, through $u^0_{;1}$, $u^0_{;2}$ and $u^0_{;3}$. However, in practice it is better to use the symmetrized system (13).

All of the above relative to the left derivative is also valid for the right derivative. However, now

$$\sigma_i = \sum_k w^j h_{ji}, \tag{11.2}$$

and $u^j_{;i} = \sigma_i \beta^j$. \hfill (12.2)

Here $w^i$ are the components of the right-hand derivative.

Noting that the left-hand sides of equations (12.1) and (12.2) are the same, combine these equations and write out in detail.

$$u^0_{;0} = \beta^0(h_{00}v^0 + h_{01}v^1 + h_{02}v^2 + h_{03}v^3) = \beta^0(h_{00}w^0 + h_{10}w^1 + h_{20}w^2 + h_{30}w^3) \quad (14.1)$$

$$u^1_{;0} = \beta^1(h_{00}v^0 + h_{01}v^1 + h_{02}v^2 + h_{03}v^3) = \beta^1(h_{00}w^0 + h_{10}w^1 + h_{20}w^2 + h_{30}w^3) \quad (14.2)$$

$$u^2_{;0} = \beta^2(h_{00}v^0 + h_{01}v^1 + h_{02}v^2 + h_{03}v^3) = \beta^2(h_{00}w^0 + h_{10}w^1 + h_{20}w^2 + h_{30}w^3) \quad (14.3)$$

$$u^3_{;0} = \beta^3(h_{00}v^0 + h_{01}v^1 + h_{02}v^2 + h_{03}v^3) = \beta^3(h_{00}w^0 + h_{10}w^1 + h_{20}w^2 + h_{30}w^3) \quad (14.4)$$

$$u^0_{;1} = \beta^0(h_{10}v^0 + h_{11}v^1 + h_{12}v^2 + h_{13}v^3) = \beta^0(h_{01}w^0 + h_{11}w^1 + h_{21}w^2 + h_{31}w^3) \quad (14.5)$$

$$u^1_{;1} = \beta^1(h_{10}v^0 + h_{11}v^1 + h_{12}v^2 + h_{13}v^3) = \beta^1(h_{01}w^0 + h_{11}w^1 + h_{21}w^2 + h_{31}w^3) \quad (14.6)$$

$$u^2_{;1} = \beta^2(h_{10}v^0 + h_{11}v^1 + h_{12}v^2 + h_{13}v^3) = \beta^2(h_{01}w^0 + h_{11}w^1 + h_{21}w^2 + h_{31}w^3) \quad (14.7)$$

$$u^3_{;1} = \beta^3(h_{10}v^0 + h_{11}v^1 + h_{12}v^2 + h_{13}v^3) = \beta^3(h_{01}w^0 + h_{11}w^1 + h_{21}w^2 + h_{31}w^3) \quad (14.8)$$

$$u^0_{;2} = \beta^0(h_{20}v^0 + h_{21}v^1 + h_{22}v^2 + h_{23}v^3) = \beta^0(h_{02}w^0 + h_{12}w^1 + h_{22}w^2 + h_{32}w^3) \quad (14.9)$$

$$u^1_{;2} = \beta^1(h_{20}v^0 + h_{21}v^1 + h_{22}v^2 + h_{23}v^3) = \beta^1(h_{02}w^0 + h_{12}w^1 + h_{22}w^2 + h_{32}w^3) \quad (14.10)$$

$$u^2_{;2} = \beta^2(h_{20}v^0 + h_{21}v^1 + h_{22}v^2 + h_{23}v^3) = \beta^2(h_{02}w^0 + h_{12}w^1 + h_{22}w^2 + h_{32}w^3) \quad (14.11)$$

$$u^3_{;2} = \beta^3(h_{20}v^0 + h_{21}v^1 + h_{22}v^2 + h_{23}v^3) = \beta^3(h_{02}w^0 + h_{12}w^1 + h_{22}w^2 + h_{32}w^3) \quad (14.12)$$

$$u^0_{;3} = \beta^0(h_{30}v^0 + h_{31}v^1 + h_{32}v^2 + h_{33}v^3) = \beta^0(h_{03}w^0 + h_{13}w^1 + h_{23}w^2 + h_{33}w^3) \quad (14.13)$$

$$u^1_{;3} = \beta^1(h_{30}v^0 + h_{31}v^1 + h_{32}v^2 + h_{33}v^3) = \beta^1(h_{03}w^0 + h_{13}w^1 + h_{23}w^2 + h_{33}w^3) \quad (14.14)$$

$$u^2_{;3} = \beta^2(h_{30}v^0 + h_{31}v^1 + h_{32}v^2 + h_{33}v^3) = \beta^2(h_{03}w^0 + h_{13}w^1 + h_{23}w^2 + h_{33}w^3) \quad (14.15)$$

$$u^3_{;3} = \beta^3(h_{30}v^0 + h_{31}v^1 + h_{32}v^2 + h_{33}v^3) = \beta^3(h_{03}w^0 + h_{13}w^1 + h_{23}w^2 + h_{33}w^3) \quad (14.16)$$

If we substitute to the left-hand side of equation (14.2) the value of $u^1_{;0}$ from the equation (13.5), reduce all three parts of the equation (14.2) by $\beta^1$ and multiply by $\beta^0$, we get the equation (14.1). Equation (14.2) is a consequence of equation (14.1). Precisely the same way we can establish that only one from each of a four equations (14.1-4), (14.5-8), (14.9-12) and (14.13-16), for example, (14.1) (14.5) (14.9) and (14.13) is independent.

But, if the equation (14.1) is multiplied by $\beta^0$, the equation (14.5) – by $\beta^1$, the equation (14.9) – by $\beta^2$, and equation (14.13) – by $\beta^3$, and add term by term, then on the left side of resultant equation will turn left side of (13.1) that is, 0.

In addition, the conditions $\sum_j h_{ij}\beta^j = \sum_j h_{ji}\beta^j = 0$ are easily verified.

Therefore, only three of equations (14.1), (14.5), (14.9) and (14.13) are linearly independent. You have to take any 3 of 4 equations and express any three of $v^i$ through the remainder, for example, through $v^0$ and through the values $u^0_{;i}$ which in these three equations are contained.

The same can be said about the values $w^i$. Resulting in the process of solving of these systems of linear equations expressions for $v^i$ and $w^i$ are too cumbersome, and we do not give them.

Thus, the components of the left and right derivatives of the function, expressed in terms of derivatives of the components are not unique, but up to arbitrary additive functions $v^0$ and $w^0$ respectively.

Let's call a function differentiable if it has left and right derivatives and these derivatives are equal.

If we require the right and left derivatives to be equal, this would mean that $v^i = w^i$ for each $0 \leq i \leq 3$.

However, in this case there are additional equations not only for the components of the derivative $v^i$, but also for the covariant derivatives of the components of a differentiable function $u^i_{;j}$.

This question requires special consideration, and we shall return to it later.

Let's take a parametrically specified isotropic curve $x^i(\tau)$ which tangent vector $\vec{\beta}$:

$$\frac{dx^i(\tau)}{d\tau} = \beta^i . \tag{15}$$

Then from equations (13.1-4) we get:

$$0 = (\sum_j u^i_{;j} \beta^j) d\tau = \sum_j u^i_{;j} dx^j = Du^i , \tag{16}$$

that is, vector $\vec{u}$ is transported in parallel along the isotropic curve $x^i(\tau)$.

If we now require this curve to be at that time the geodesic one, we get, by substituting (15) into the geodesic equation:

$$\frac{d^2 x^i}{d\tau^2} + \sum_{k,l} \Gamma^i_{kl} \frac{dx^k}{d\tau} \frac{dx^l}{d\tau} = \frac{d\beta^i}{d\tau} + \sum_{k,l} \Gamma^i_{kl} \beta^k \beta^l =$$

$$= \sum_l \frac{\partial \beta^i}{\partial x^l} \frac{dx^l}{d\tau} + \sum_{k,l} \Gamma^i_{kl} \beta^k \beta^l = \sum_l (\frac{\partial \beta^i}{\partial x^l} + \sum_k \Gamma^i_{kl} \beta^k) \beta^l = . \tag{17}$$

$$= \sum_l \beta^i_{;l} \beta^l = 0$$

There are the null geodesic equations.

Now, compare this expression with the first 4 analogues of the Cauchy-Riemann equations (13):
$\sum_j u^i_{;j} \beta^j = 0$.

So we have established that differentiable from the right (or left) vector-function $\vec{u}$ is transported in parallel along the isotropic line which tangent vector is determined by vector $\vec{\beta}$.

The angle between the vectors $\vec{u}$ and $\vec{\beta}$, of course, always remains constant. But then it is not clear what for the initial isotropic vector $\vec{\alpha}$ in the tangent space having the appointed (but arbitrary) for given algebra components $\alpha^i$ was introduced. And what could mean a corresponding isotropic direction of the vector $\vec{\beta}$?

It remains to make the final step. We assume that the vector $\vec{u}$ *is* vector $\vec{\beta}$:

$\vec{u} \equiv \vec{\beta}$. (18)

Then the first 4 analogues of the Cauchy-Riemann equations are the equations is an isotropic geodesic.

And equations (13) are converted into the equations:

$$\sum_j \beta^i_{;j} \beta^j = 0 , \quad 0 \leq i \leq 3$$

$$\beta^0 \beta^i_{;j} = \beta^i \beta^0_{;j}, \quad 1 \leq i \leq 3, \quad 0 \leq j \leq 3 \tag{19}$$

And equations (13a) are converted into the equations:

$$\sum_j \beta^i_{;j} \beta^j = 0 , \quad 0 \leq i \leq 3 \tag{19a}$$

$$\beta^j_{;i} \beta^k = \beta^k_{;i} \beta^j, \quad 0 \leq i,j,k \leq 3$$

Equations (13.1-16) (note that they are just in detail painted equations (13)) are converted into the equations

$$\beta^0_{;0}\beta^0 + \beta^0_{;1}\beta^1 + \beta^0_{;2}\beta^2 + \beta^0_{;3}\beta^3 = 0 \tag{19.1}$$

$$\beta^1_{;0}\beta^0 + \beta^1_{;1}\beta^1 + \beta^1_{;2}\beta^2 + \beta^1_{;3}\beta^3 = 0 \tag{19.2}$$

$$\beta^2_{;0}\beta^0 + \beta^2_{;1}\beta^1 + \beta^2_{;2}\beta^2 + \beta^2_{;3}\beta^3 = 0 \tag{19.3}$$

$$\beta^3_{;0}\beta^0 + \beta^3_{;1}\beta^1 + \beta^3_{;2}\beta^2 + \beta^3_{;3}\beta^3 = 0 \tag{19.4}$$

$$\beta^1_{;0} = \frac{\beta^1}{\beta^0} \beta^0_{;0} \tag{19.5}$$

$$\beta^2_{;0} = \frac{\beta^2}{\beta^0} \beta^0_{;0} \tag{19.6}$$

$$\beta^3_{;0} = \frac{\beta^3}{\beta^0} \beta^0_{;0} \tag{19.7}$$

$$\beta^1_{;1} = \frac{\beta^1}{\beta^0} \beta^0_{;1} \tag{19.8}$$

$$\beta^2_{;1} = \frac{\beta^2}{\beta^0} \beta^0_{;1} \tag{19.9}$$

$$\beta^3_{;1} = \frac{\beta^3}{\beta^0} \beta^0_{;1} \tag{19.10}$$

$$\beta^1_{;2} = \frac{\beta^1}{\beta^0} \beta^0_{;2} \tag{19.11}$$

$$\beta^2_{;2} = \frac{\beta^2}{\beta^0} \beta^0_{;2} \tag{19.12}$$

$$\beta^3_{;2} = \frac{\beta^3}{\beta^0} \beta^0_{;2} \tag{19.13}$$

$$\beta^1_{;3} = \frac{\beta^1}{\beta^0} \beta^0_{;3} \tag{19.14}$$

$$\beta^2_{;3} = \frac{\beta^2}{\beta^0} \beta^0_{;3} \tag{19.15}$$

$$\beta^3_{;3} = \frac{\beta^3}{\beta^0} \beta^0_{;3} \tag{19.16}$$

Now recall the definition of the vector $\vec{\beta}$:

$$\beta^i = \sum_j \alpha^j b^i_j \tag{7}$$

Here $\vec{\alpha} = (1, \alpha^1(\vec{r}), \alpha^2(\vec{r}), \alpha^3(\vec{r}))$ and matrix $b^i_j$ is a function of the metric tensor only, which in it's turn depends on the coordinates.

Thus, 13 independent (as we found earlier) equations (13) and, consequently, 13 independent equations (19) is the system of quasi-linear differential equations of the 1st order for 10 independent components of the metric tensor and 3 unknown components $\alpha^1$, $\alpha^2$ и $\alpha^3$ for a field of isotropic directions.

However, to the equations (19) we must add the constraint equation
$$(\alpha^1)^2 + (\alpha^2)^2 + (\alpha^3)^2 = 1. \qquad (2)$$
Thus, the system (2), (19) (or (2), (19a)) is overdetermined.
In general, the system looks extremely difficult to be solved.
However, we can solve this system for the metric tensor and the vector $\vec{\beta}$ immediately without having deal with it for the vector $\vec{\alpha}$.
The constraint (2) turns into the equation
$$\sum_{i,j} g_{ij} \beta^i \beta^j = 0. \qquad (20)$$
And the condition $\alpha^0 = 1$ turns into the condition
$$1 = \alpha^0 = \sum_j \beta^j a^0_j. \qquad (21)$$
Thus, the conditions (20) and (21) reduce the number of independent components of the vector $\vec{\beta}$ to two.

But in view of quadratic terms of the condition (20) it does not seem reasonable to express two any components of the vector $\vec{\beta}$ through the other two with the help of the conditions (20) and (21), and subsequent substitution of these values into equations (19) or (19a), as it brakes the symmetric record of them.

Now let's come back to the question of what will happen if the left and right derivatives are equal, that is $w^i = v^i$.

In the future, once we have decided to work with the values $\beta^i$ only (and not with $\alpha^i$), then the values $h_{ij} = \sum_{m,n} a^m_i a^n_j \kappa_{ij}$, the definition of which includes values $\kappa_{ij}$, which in turn, contain $\alpha^i$, should be recalculated taking into account the expression
$$\alpha^i = \sum_j \beta^j a^i_j.$$
As noted earlier, the conditions
$$\sum_j h_{ij} \beta^j = \sum_j h_{ji} \beta^j = 0 \qquad (22)$$
can be easily verified.
It is also easy to verify the conditions
$$h_{0i} = h_{i0}. \qquad (23)$$
As it will be seen below, these conditions are essential.
In particular, for example, when $i = 1$ conditions (22) give
$$h_{10}\beta^0 + h_{11}\beta^1 + h_{12}\beta^2 + h_{13}\beta^3 = h_{01}\beta^0 + h_{11}\beta^1 + h_{21}\beta^2 + h_{31}\beta^3.$$
Or, according to the conditions (23):
$$(h_{12} - h_{21})\beta^2 = -(h_{13} - h_{31})\beta^3.$$
Now let's consider the right-hand equality of triple equation (14.5), of course, not forgetting that now $w^i = v^i$ and $\vec{u} \equiv \vec{\beta}$:
$$\beta^0_{;1} = \beta^0(h_{10}v^0 + h_{11}v^1 + h_{12}v^2 + h_{13}v^3) = \beta^0(h_{01}v^0 + h_{11}v^1 + h_{21}v^2 + h_{31}v^3),$$

or, using (23) and considering only the right-hand equality of this triple equation,
$(h_{12} - h_{21})v^2 = -(h_{13} - h_{31})v^3$.
Dividing the left side of the penultimate equation by the left side of the latter, and the right-hand side - by the right, respectively, and equating the obtained expressions, we find
$$\frac{\beta^2}{v^2} = \frac{\beta^3}{v^3}.$$
Acting similarly with the right-hand equalities of triple equations (14.9) and (14.13), we find the relations $\frac{\beta^1}{v^1} = \frac{\beta^3}{v^3}$ и $\frac{\beta^1}{v^1} = \frac{\beta^2}{v^2}$.

From the last three ratios we can express two of three components of derivative through the third, for example

$$v^2 = \frac{\beta^2}{\beta^1}v^1, \qquad (24)$$

$$v^3 = \frac{\beta^3}{\beta^1}v^1. \qquad (25)$$

Consider in more detail the consequences of the relations (24) and (25).
We take the equations (14.1) and substitute into them the values of $v^2$ and $v^3$ from (24) and (25), again remembering that now $w^i = v^i$ and $\vec{u} \equiv \vec{\beta}$.

$$\beta^0_{;0} = \beta^0(h_{00}v^0 + h_{01}v^1 + h_{02}v^2 + h_{03}v^3) = \beta^0(h_{00}v^0 + h_{10}v^1 + h_{20}v^2 + h_{30}v^3) =$$
$$= \beta^0(h_{00}v^0 + h_{01}v^1 + h_{02}\frac{\beta^2}{\beta^1}v^1 + h_{03}\frac{\beta^3}{\beta^1}v^1) = \beta^0(h_{00}v^0 + h_{10}v^1 + h_{20}\frac{\beta^2}{\beta^1}v^1 + h_{30}\frac{\beta^3}{\beta^1}v^1) =$$
$$= \beta^0(h_{00}v^0 + \frac{h_{01}\beta^1 + h_{02}\beta^2 + h_{03}\beta^3}{\beta^1}v^1) = \beta^0(h_{00}v^0 + \frac{h_{10}\beta^1 + h_{20}\beta^2 + h_{30}\beta^3}{\beta^1}v^1) =$$
$$= \beta^0(h_{00}v^0 - \frac{h_{00}\beta^0}{\beta^1}v^1) \equiv \beta^0(h_{00}v^0 - \frac{h_{00}\beta^0}{\beta^1}v^1) = \beta^0 h_{00}(v^0 - \frac{\beta^0}{\beta^1}v^1) \equiv \beta^0 h_{00}(v^0 - \frac{\beta^0}{\beta^1}v^1)$$

The equality in the penultimate line is carried out by the conditions (22).
Similarly, we can transform the equations (14.2):
$$\beta^1_{;0} = \beta^1(h_{00}v^0 + h_{01}v^1 + h_{02}v^2 + h_{03}v^3) = \beta^1(h_{00}v^0 + h_{10}v^1 + h_{20}v^2 + h_{30}v^3) =$$
$$= \beta^1 h_{00}(v^0 - \frac{\beta^0}{\beta^1}v^1) \equiv \beta^1 h_{00}(v^0 - \frac{\beta^0}{\beta^1}v^1)$$

And, in general, if we introduce the designation
$$S = v^0 - \frac{\beta^0}{\beta^1}v^1, \qquad (26)$$

we obtain while transforming equations (14.1-16):
$$\beta^i_{;j} = \beta^i h_{j0} S \equiv \beta^i h_{0j} S. \qquad (27)$$

The last identity takes place because of (23).
We write the equations (27) in more detail.

$\beta^0_{;0} = \beta^0 h_{00} S$ (27.1)    $\beta^1_{;0} = \beta^1 h_{00} S$ (27.2)    $\beta^2_{;0} = \beta^2 h_{00} S$ (27.3)    $\beta^3_{;0} = \beta^3 h_{00} S$ (27.4)

$\beta^0_{;1} = \beta^0 h_{10} S$ (27.5)    $\beta^1_{;1} = \beta^1 h_{10} S$ (27.6)    $\beta^2_{;1} = \beta^2 h_{10} S$ (27.7)    $\beta^3_{;1} = \beta^3 h_{10} S$ (27.8)

$\beta^0_{;2} = \beta^0 h_{20} S$ (27.9)    $\beta^1_{;2} = \beta^1 h_{20} S$ (27.10)    $\beta^2_{;2} = \beta^2 h_{20} S$ (27.11)    $\beta^3_{;2} = \beta^3 h_{20} S$ (27.12)

$\beta^0_{;3} = \beta^0 h_{30} S$ (27.13)    $\beta^1_{;3} = \beta^1 h_{30} S$ (27.14)    $\beta^2_{;3} = \beta^2 h_{30} S$ (27.15)    $\beta^3_{;3} = \beta^3 h_{30} S$ (27.16)

We note that if the equation (27.1) is multiplied by $\beta^0$, (27.5) - by $\beta^1$, (27.9) – by $\beta^2$, (27.13) – by $\beta^3$ and add term by term, then the left side of the obtained equation turns to the left-hand side of the Cauchy-Riemann condition (19.1), and the right side because of (22) becomes zero.

Exactly the same results can be obtained for (19.2-4), if we perform the same actions with the other three columns of (27.1-16).

It follows that for a differentiable function the Cauchy-Riemann equations (19.1-4) or the first four equations (19) or (19a), which are at the same time null geodetic equations are performed *identically*:

$$\sum_j \beta^i_{;j} \beta^j \equiv 0 . \tag{28}$$

If we now make all the possible divisions of the right and left sides of the equations (27.1), (27.5), (27.9) and (27.13), we obtain:

$$\frac{\beta^0_{;0}}{\beta^0_{;1}} = \frac{h_{00}}{h_{10}} \qquad \frac{\beta^0_{;0}}{\beta^0_{;2}} = \frac{h_{00}}{h_{20}} \qquad \frac{\beta^0_{;0}}{\beta^0_{;3}} = \frac{h_{00}}{h_{30}}$$

$$\frac{\beta^0_{;1}}{\beta^0_{;0}} = \frac{h_{10}}{h_{00}} \qquad \frac{\beta^0_{;1}}{\beta^0_{;2}} = \frac{h_{10}}{h_{20}} \qquad \frac{\beta^0_{;1}}{\beta^0_{;3}} = \frac{h_{10}}{h_{30}}$$

$$\frac{\beta^0_{;2}}{\beta^0_{;0}} = \frac{h_{20}}{h_{00}} \qquad \frac{\beta^0_{;2}}{\beta^0_{;1}} = \frac{h_{20}}{h_{10}} \qquad \frac{\beta^0_{;2}}{\beta^0_{;3}} = \frac{h_{20}}{h_{30}}$$

$$\frac{\beta^0_{;3}}{\beta^0_{;0}} = \frac{h_{30}}{h_{00}} \qquad \frac{\beta^0_{;3}}{\beta^0_{;1}} = \frac{h_{30}}{h_{10}} \qquad \frac{\beta^0_{;3}}{\beta^0_{;2}} = \frac{h_{30}}{h_{20}}$$

or

$$\beta^0_{;0} h_{10} = \beta^0_{;1} h_{00} \qquad \beta^0_{;0} h_{20} = \beta^0_{;2} h_{00} \qquad \beta^0_{;0} h_{30} = \beta^0_{;3} h_{00}$$
$$\beta^0_{;1} h_{00} = \beta^0_{;0} h_{10} \qquad \beta^0_{;1} h_{20} = \beta^0_{;2} h_{10} \qquad \beta^0_{;1} h_{30} = \beta^0_{;3} h_{10}$$
$$\beta^0_{;2} h_{00} = \beta^0_{;0} h_{20} \qquad \beta^0_{;2} h_{10} = \beta^0_{;1} h_{20} \qquad \beta^0_{;2} h_{30} = \beta^0_{;3} h_{20}$$
$$\beta^0_{;3} h_{00} = \beta^0_{;0} h_{30} \qquad \beta^0_{;3} h_{10} = \beta^0_{;1} h_{30} \qquad \beta^0_{;3} h_{20} = \beta^0_{;2} h_{30} .$$

If we add here the trivial relations $\beta^0_{;i} h_{i0} = \beta^0_{;i} h_{i0}$, we get:

$$\beta^0_{;i} h_{j0} = \beta^0_{;j} h_{i0} .$$

The same relations can be derived for the other three columns of (27.1-16). Then we get:

$$\beta^k_{;i} h_{j0} = \beta^k_{;j} h_{i0} \qquad i, j, k = 0,...,3 \tag{29}$$

If we now, for example, multiply the equation (27.1) by $h_{00}$, (27.2) – by $h_{01}$, (27.3) – by $h_{02}$, (27.4) – by $h_{03}$ and add them term by term, we obtain

$$\beta^0_{;0} h_{00} + \beta^1_{;0} h_{01} + \beta^2_{;0} h_{02} + \beta^3_{;0} h_{03} = h_{00} S(\beta^0 h_{00} + \beta^1 h_{01} + \beta^2 h_{02} + \beta^3 h_{03}) = 0 . \tag{30}$$

The last equality takes place because of (22).

If we continue with this kind of reasoning, it would seem that there are new equations

$$\sum_m \beta^m_{;j} h_{im} = 0 \text{ и } \sum_m \beta^m_{;j} h_{im} = 0 . \tag{31}$$

But if we substitute the values of $\beta^1_{;0}$, $\beta^2_{;0}$ and $\beta^3_{;0}$ from the equations (19.5-7) into the left-hand side of equation (30), we obtain

$$\beta^0_{;0} h_{00} + \beta^1_{;0} h_{01} + \beta^2_{;0} h_{02} + \beta^3_{;0} h_{03} = \beta^0_{;0} h_{00} + \frac{\beta^1}{\beta^0} \beta^0_{;0} h_{01} + \frac{\beta^2}{\beta^0} \beta^0_{;0} h_{02} + \frac{\beta^3}{\beta^0} \beta^0_{;0} h_{03} =$$

$$= \frac{1}{\beta^0} \beta^0_{;0} (\beta^0 h_{00} + \beta^1 h_{01} + \beta^2 h_{02} + \beta^3 h_{03}) = 0$$

If we continue to act the same way, we find that equations (31) become identities.

If we consider the system (27.1-16) as a system of linear equations for the 16 unknown $\beta^i_{;j}$ and the unknown $S$, then 15 of 16 values of $\beta^i_{;j}$ can be expressed through any one of them, $\beta^0_{;0}$ for example.

Then we obtain:

$$\beta^1_{;0} = \frac{\beta^1}{\beta^0} \beta^0_{;0} \tag{32.1}$$

$$\beta^2_{;0} = \frac{\beta^2}{\beta^0} \beta^0_{;0} \tag{32.2}$$

$$\beta^3_{;0} = \frac{\beta^3}{\beta^0} \beta^0_{;0} \tag{32.3}$$

$$\beta^0_{;1} = \frac{h_{10}}{h_{00}} \beta^0_{;0} \tag{32.4}$$

$$\beta^1_{;1} = \frac{\beta^1 h_{10}}{\beta^0 h_{00}} \beta^0_{;0} \tag{32.5}$$

$$\beta^2_{;1} = \frac{\beta^2 h_{10}}{\beta^0 h_{00}} \beta^0_{;0} \tag{32.6}$$

$$\beta^3_{;1} = \frac{\beta^3 h_{10}}{\beta^0 h_{00}} \beta^0_{;0} \tag{32.7}$$

$$\beta^0_{;2} = \frac{h_{20}}{h_{00}} \beta^0_{;0} \tag{32.8}$$

$$\beta^1_{;2} = \frac{\beta^1 h_{20}}{\beta^0 h_{00}} \beta^0_{;0} \tag{32.9}$$

$$\beta^2_{;2} = \frac{\beta^2 h_{20}}{\beta^0 h_{00}} \beta^0_{;0} \tag{32.10}$$

$$\beta^3_{;2} = \frac{\beta^3 h_{20}}{\beta^0 h_{00}} \beta^0_{;0} \tag{32.11}$$

$$\beta^0_{;3} = \frac{h_{30}}{h_{00}} \beta^0_{;0} \tag{32.12}$$

$$\beta^1_{;3} = \frac{\beta^1 h_{30}}{\beta^0 h_{00}} \beta^0_{;0} \tag{32.13}$$

$$\beta^2_{;3} = \frac{\beta^2 h_{30}}{\beta^0 h_{00}} \beta^0_{;0} \tag{32.14}$$

$$\beta^3_{;3} = \frac{\beta^3 h_{30}}{\beta^0 h_{00}} \beta^0_{;0} \tag{32.15}$$

Or in general, including trivial relations:

$$\beta^i_{;j} = \frac{\beta^i h_{j0}}{\beta^0 h_{00}} \beta^0_{;0} \qquad i, j = 0,...,3. \tag{32}$$

In this case (and from (26)), we have:

$$S = v^0 - \frac{\beta^0}{\beta^1} v^1 = \frac{\beta^0_{;0}}{\beta^0 h_{00}}.$$

From this relation, together with (24) and (25) we finally find the connection between the components of the derivative of a differentiable function:

$$v^1 = \frac{\beta^1}{\beta^0}v^0 - \frac{\beta^1}{(\beta^0)^2 h_{00}}u^0_{;0}, \tag{33.1}$$

$$v^2 = \frac{\beta^2}{\beta^0}v^0 - \frac{\beta^2}{(\beta^0)^2 h_{00}}u^0_{;0}, \tag{33.2}$$

$$v^3 = \frac{\beta^3}{\beta^0}v^0 - \frac{\beta^3}{(\beta^0)^2 h_{00}}u^0_{;0}. \tag{33.3}$$

We summarize the above.

**1)**. Covariant derivatives of the components of differentiable from the right and/or from the left vector function in this pseudo-Riemannian space (in the case when this vector function *determines the direction* of zero geodetic) obeys the system of equations

$$\sum_j \beta^i_{;j}\beta^j = 0, \quad 0 \leq i \leq 3$$
$$\beta^0 \beta^i_{;j} = \beta^i \beta^0_{;j}, \quad 1 \leq i \leq 3, \quad 0 \leq j \leq 3 \tag{19}$$

or

$$\sum_j \beta^i_{;j}\beta^j = 0, \quad 0 \leq i \leq 3$$
$$\beta^j_{;i}\beta^k = \beta^k_{;i}\beta^j, \quad 0 \leq i,j,k \leq 3 \tag{19a}$$

and

$$\sum_{i,j} g_{ij}\beta^i \beta^j = 0, \tag{20}$$

$$\sum_j \beta^j a^0_j = 1. \tag{21}$$

Thus differentiable from the right and/or from the left vector function has two independent components.

The components of the derivative of a differentiable from the right and/or from the left vector function $v^i$ and $w^i$ are evaluated to an approximation of additive function.

It is necessary to take any three of the following equations

$$\beta^0_{;0} = \beta^0(h_{00}v^0 + h_{01}v^1 + h_{02}v^2 + h_{03}v^3) = \beta^0(h_{00}w^0 + h_{10}w^1 + h_{20}w^2 + h_{30}w^3)$$
$$\beta^0_{;1} = \beta^0(h_{10}v^0 + h_{11}v^1 + h_{12}v^2 + h_{13}v^3) = \beta^0(h_{01}w^0 + h_{11}w^1 + h_{21}w^2 + h_{31}w^3)$$
$$\beta^0_{;2} = \beta^0(h_{20}v^0 + h_{21}v^1 + h_{22}v^2 + h_{23}v^3) = \beta^0(h_{02}w^0 + h_{12}w^1 + h_{22}w^2 + h_{32}w^3)$$
$$\beta^0_{;3} = \beta^0(h_{30}v^0 + h_{31}v^1 + h_{32}v^2 + h_{33}v^3) = \beta^0(h_{03}w^0 + h_{13}w^1 + h_{23}w^2 + h_{33}w^3)$$

and to express any 3 of the 4 values of $v^i$ through the remainder, for example, through $v^0$ and through quantities $u^0_{;i}$ containing in these three equations. The same can be said about the values of $w^i$.

**2).** Covariant derivatives of the components of *differentiable* vector function in this pseudo-Riemannian space (in the case when the vector function *determines the direction* of zero geodesic) obeys the system of equations

$$\beta^0 \beta^i_{;j} = \beta^i \beta^0_{;j}, \quad 1 \leq i \leq 3, \quad 0 \leq j \leq 3, \tag{19c}$$

$$\beta^k_{;i}h_{j0} = \beta^k_{;j}h_{i0} \quad i,j,k = 0,...,3, \tag{29}$$

$$\sum_{i,j} g_{ij}\beta^i \beta^j = 0, \tag{20}$$

$$\sum_j \beta^j a^0_j = 1. \tag{21}$$

or recorded in less symmetrized form equivalent system

$$\beta^i_{;j} = \frac{\beta^i h_{j0}}{\beta^0 h_{00}} \beta^0_{;0} \qquad i,j = 0,...,3, \tag{32}$$

$$\sum_{i,j} g_{ij} \beta^i \beta^j = 0, \tag{20}$$

$$\sum_j \beta^j a^0_j = 1. \tag{21}$$

The components of the derivative of a differentiable vector function has two independent components (as the initial vector function), which can be seen from the following relations

$$v^1 = \frac{\beta^1}{\beta^0} v^0 - \frac{\beta^1}{(\beta^0)^2 h_{00}} u^0_{;0}, \tag{33.1}$$

$$v^2 = \frac{\beta^2}{\beta^0} v^0 - \frac{\beta^2}{(\beta^0)^2 h_{00}} u^0_{;0}, \tag{33.2}$$

$$v^3 = \frac{\beta^3}{\beta^0} v^0 - \frac{\beta^3}{(\beta^0)^2 h_{00}} u^0_{;0}. \tag{33.3}$$

Let us try to solve the system (32), (20), (21) in the simplest cases.
In the case of Minkowski geometry the solution can be found instantly:
$\vec{\beta} = (1, \beta^1, \beta^2, \beta^3)$, where $\beta^i = const$ and $(\beta^1)^2 + (\beta^2)^2 + (\beta^3)^2 = 1$, i.e. the cone of light that, in general, is not surprising.
Let's consider the case of the spherical metric $g = diag(L^2(t,r), -H^2(t,r), -r^2 \sin^2(\vartheta), -r^2)$.
Independent coordinates are indicated in the following order: $x^1 = t$, $x^2 = r$, $x^3 = \varphi$, $x^4 = \vartheta$.
Equation (21) immediately gives the dependence

$$\beta^1(t,r,\varphi,\vartheta) = \frac{1}{L(t,r)}.$$

If we substitute this value of $\beta^1(t,r,\varphi,\vartheta)$ into the equation (32.8), we obtain (to an approximation of a factor that can not be zero)
$L_r(t,r) = 0$.
From here we obtain
$L(t,r) = L^{(1)}(t)$.
Now, of course,

$$\beta^1(t,r,\varphi,\vartheta) = \frac{1}{L^{(1)}(t)}.$$

Here $L^{(1)}$ is a symbol of a new variable.
We substitute these new values of $L(t,r)$ and $\beta^1(t,r,\varphi,\vartheta)$ into all equations. In the future all the time we act in this way.
Equation (32.4) gives (again to an approximation of nonzero factor)
$H_t(t,r) = 0$,
from where
$H(t,r) = H^{(1)}(r)$,
where $H^{(1)}(r)$ is again a new variable.
Equation (32.2) gives
$\beta^3_t(t,r,\varphi,\vartheta) = 0$,
from which we immediately obtain
$\beta^3(t,r,\varphi,\vartheta) = \beta^{3(1)}(r,\varphi,\vartheta)$.
Absolutely in the same way equation (32.3) gives

$\beta^4(t,r,\varphi,\vartheta) = \beta^{4(1)}(r,\varphi,\vartheta)$.

Next, in the same way equation (32.1) gives

$\beta^2(t,r,\varphi,\vartheta) = \beta^{2(1)}(r,\varphi,\vartheta)$.

Equation (32.6) is converted to

$$\frac{\beta^{3(1)}(r,\varphi,\vartheta)}{r} + \beta_r^{3(1)}(r,\varphi,\vartheta) = 0,$$

from which we obtain

$$\beta^{3(1)}(r,\varphi,\vartheta) = \frac{\beta^{3(2)}(\varphi,\vartheta)}{r}.$$

Quite similar the equation (32.7) gives

$$\beta^{4(1)}(r,\varphi,\vartheta) = \frac{\beta^{4(2)}(\varphi,\vartheta)}{r}.$$

Further substitutions into the system yields the following.
Equation (32.5) is converted into equation

$$\beta^{2(1)}(r,\varphi,\vartheta)H_r(r) + \beta_r^{2(1)}(r,\varphi,\vartheta)H(r) = 0,$$

from which we obtain

$$\beta^{2(1)}(r,\varphi,\vartheta) = \frac{\beta^{2(2)}(\varphi,\vartheta)}{H(r)}.$$

Equation (32.14) takes the form

$\cot\vartheta\,\beta^{3(2)}(\varphi,\vartheta) + \beta_\vartheta^{3(2)}(\varphi,\vartheta) = 0$,

from which we get

$\beta^{3(2)}(\varphi,\vartheta) = \csc\vartheta\,\beta^{3(3)}(\varphi)$.

Next, equation (32.15) takes the form

$\beta^{2(2)}(\varphi,\vartheta) + H^{(1)}(r)\beta_\vartheta^{4(2)}(\varphi,\vartheta) = 0$,

from where, for obvious reasons, follows

$H^{(1)}(r) = H^{(2)} = const$.

Further equation (32.15) takes the form

$\beta^{2(2)}(\varphi,\vartheta) + H^{(2)}\beta_\vartheta^{4(2)}(\varphi,\vartheta) = 0$.

And equation (32.13) takes the form

$-\beta^{4(2)}(\varphi,\vartheta) + H^{(2)}\beta_\vartheta^{2(2)}(\varphi,\vartheta) = 0$.

Now, if the first of them we differentiate by $\vartheta$ and multiply by $H^{(2)}$, and then subtract the resulting equation from the second, we obtain the equation

$\beta^{4(2)}(\varphi,\vartheta) + (H^{(2)})^2\beta_{\vartheta\vartheta}^{4(2)}(\varphi,\vartheta) = 0$,

which is easily integrated:

$$\beta^{4(2)}(\varphi,\vartheta) = \cos\frac{\vartheta}{H^{(2)}}\gamma(\varphi) + \sin\frac{\vartheta}{H^{(2)}}\delta(\varphi),$$

where $\gamma(\varphi)$ and $\delta(\varphi)$ are arbitrary functions.
After that equation (32.13) takes the form

$$-\cos\frac{\vartheta}{H^{(2)}}\gamma(\varphi) - \sin\frac{\vartheta}{H^{(2)}}\delta(\varphi) + H^{(2)}\beta_\vartheta^{2(2)}(\varphi,\vartheta) = 0.$$

This equation is also easily integrated:

$$\beta^{2(2)}(\varphi,\vartheta) = \sin\frac{\vartheta}{H^{(2)}}\gamma(\varphi) - \cos\frac{\vartheta}{H^{(2)}}\delta(\varphi) + \varepsilon(\varphi),$$

where $\varepsilon(\varphi)$ is an arbitrary function.
As a result, equation (32.15) takes the form (we write now, for clarity, all the factors that can not be zero):

$$\frac{a(t,r,\varphi,\vartheta)L^{(1)}(t)\varepsilon(\varphi)}{H^{(2)}r}=0.$$

Obviously, there is the only possibility that the task does not become a trivial:
$\varepsilon(\varphi)=0$.

Finally, we write the resulting equation (32.9) and (32.11), again to an approximation of nonzero factors:

$$-\sin\vartheta\,\beta^{3(3)}(\varphi)+H^{(2)}\sin\frac{\vartheta}{H^{(2)}}\gamma_\varphi(\varphi)-H^{(2)}\cos\frac{\vartheta}{H^{(2)}}\delta_\varphi(\varphi)=0,$$

$$-\cos\vartheta\,\beta^{3(3)}(\varphi)+\cos\frac{\vartheta}{H^{(2)}}\gamma_\varphi(\varphi)+\sin\frac{\vartheta}{H^{(2)}}\delta_\varphi(\varphi)=0.$$

If now we the first one multiply by $\cos\vartheta$, and the second by $\sin\vartheta$, and if we subtract from the first obtained equation the second one, we obtain

$$(-\cos\frac{\vartheta}{H^{(2)}}\sin\vartheta+H^{(2)}\cos\vartheta\sin\frac{\vartheta}{H^{(2)}})\gamma_\varphi(\varphi)-(\sin\frac{\vartheta}{H^{(2)}}\sin\vartheta+H^{(2)}\cos\vartheta\cos\frac{\vartheta}{H^{(2)}})\delta_\varphi(\varphi)=0$$

or

$$\frac{-\cos\frac{\vartheta}{H^{(2)}}\sin\vartheta+H^{(2)}\cos\vartheta\sin\frac{\vartheta}{H^{(2)}}}{\sin\frac{\vartheta}{H^{(2)}}\sin\vartheta+H^{(2)}\cos\vartheta\cos\frac{\vartheta}{H^{(2)}}}=\frac{\delta_\varphi(\varphi)}{\gamma_\varphi(\varphi)}.$$

In the last equation we see that the left side is a function of $\vartheta$ only, and the right – of $\varphi$ only. Consequently, the left and right sides must be equal to the same constant. It is clear that this fact can not take place.

So, spherically symmetric solution does not exist.

Now we consider a cylindrically symmetric metric
$g=diag(L^2(t,r),-H^2(t,r),-K^2(t,r),-P^2(t,r))$.
Independent coordinates are indicated in the following order: $x^1=t$, $x^2=r$, $x^3=\varphi$, $x^4=z$.
Just as in the case of a spherically symmetry, the equation (21) immediately gives dependence
$$\beta^1(t,r,\varphi,\vartheta)=\frac{1}{L(t,r)}.$$
Equation (32.4) is converted into equation
$H(t,r)\beta^2(t,r,\varphi,z)L_r(t,r)+H_t(t,r)=0$,
from which it is clear that
$\beta^2(t,r,\varphi,z)=\beta^{2(1)}(t,r)$.
Equation (32.2) now looks like this:
$-K(t,r)\beta^{2(1)}(t,r)\beta^3(t,r,\varphi,z)L_r(t,r)+\beta^3(t,r,\varphi,z)K_t(t,r)+K(t,r)\beta^3_t(t,r,\varphi,z)=0$.
If the left-hand side of this equation we divide by $\beta^3(t,r,\varphi,z)$ and by $K(t,r)$ we get
$$-\beta^{2(1)}(t,r)L_r(t,r)+\frac{K_t(t,r)}{K(t,r)}=-\frac{\beta^3_t(t,r,\varphi,z)}{\beta^3(t,r,\varphi,z)},$$
from where it is clear that the right side must be a function of $t$ and $r$ only. We introduce a new function $\beta^{3(1)}(t,r)$, and solve the equation
$$\frac{\beta^3_t(t,r,\varphi,z)}{\beta^3(t,r,\varphi,z)}=\beta^{3(1)}(t,r).$$
The solution of this equation is
$\beta^3(t,r,\varphi,z)=e^{\beta^{3(1)}(t,r)}\beta^{3(2)}(r,\varphi,z)$,
where $\beta^{3(2)}(r,\varphi,z)$ is an arbitrary function.

Proceeding in a similar manner the equation (32.6), it appears that
$$\beta^{3(2)}(r,\varphi,z) = e^{\beta^{3(3)}(t,r)}\beta^{3(4)}(\varphi,z),$$
from where
$$\beta^3(t,r,\varphi,z) = e^{\beta^{3(1)}(t,r)+\beta^{3(3)}(t,r)}\beta^{3(4)}(\varphi,z),$$
where $\beta^{3(4)}(\varphi,z)$ is an arbitrary function.
Renaming variables, we finally find
$$\beta^3(t,r,\varphi,z) = e^{\beta^{3(1)}(t,r)}\beta^{3(2)}(\varphi,z).$$
By doing in precisely the same manner with equations (32.7) and (32.3), we find out that
$$\beta^4(t,r,\varphi,z) = e^{\beta^{4(1)}(t,r)}\beta^{4(2)}(\varphi,z).$$
Then equations (32.4), (32.8) and (32.12) become respectively
$$H(t,r)\beta^{2(1)}(t,r)L_r(t,r) + H_t(t,r) = 0,$$
$$K(t,r)\beta^{2(1)}(t,r)L_r(t,r) + K_t(t,r) = 0,$$
$$P(t,r)\beta^{2(1)}(t,r)L_r(t,r) + P_t(t,r) = 0.$$
Multiplying the first of these equations by $K(t,r)$, the second by $-H(t,r)$, and adding the resulting equations, we obtain
$$K(t,r)H_t(t,r) - H(t,r)K_t(t,r) = 0.$$
Integrating this equation, we get
$$K(t,r) = H(t,r)K^{(1)}(r).$$
Quite similarly, we obtain
$$P(t,r) = H(t,r)P^{(1)}(r).$$
Here $K^{(1)}(r)$ and $P^{(1)}(r)$ are arbitrary functions.
Now the equation (32.9) and (32.13) look like this:
$$-H(t,r)L(t,r)K_r^{(1)}(r) - K^{(1)}(r)L(t,r)H_r(t,r) + (H(t,r))^3 K^{(1)}(r)(\beta^{2(1)}(t,r))^2 L_r(t,r) = 0,$$
$$-H(t,r)L(t,r)P_r^{(1)}(r) - P^{(1)}(r)L(t,r)H_r(t,r) + (H(t,r))^3 P^{(1)}(r)(\beta^{2(1)}(t,r))^2 L_r(t,r) = 0.$$
Multiplying the first of them by $-P^{(1)}(r)$, the second by $K^{(1)}(r)$, and adding the resulting equations, we have to an approximation of a non-zero factor
$$P^{(1)}(r)K_r^{(1)}(r) - K^{(1)}(r)P_r^{(1)}(r) = 0.$$
From where once we get
$$P^{(1)}(r) = P^{(2)}K^{(1)}(r),$$
where $P^{(2)} = const$.
After that, we make the substitution
$$H(t,r) = e^{H^{(1)}(t,r)}.$$
Equations (32.12) and (32.2) take the form:
$$\beta^{2(1)}(t,r)L_r(t,r) + H_t^{(1)}(t,r) = 0,$$
$$\beta^{2(1)}(t,r)L_r(t,r) - H_t^{(1)}(t,r) - \beta_t^{3(1)}(t,r) = 0.$$
Subtracting the second equation from the first and integrating the resulting equation, we find
$$\beta^{3(1)}(t,r) = -2 H^{(1)}(t,r) + \beta^{3(3)}(r),$$
where $\beta^{3(3)}(r)$ is an arbitrary function.
By doing the same with a pair of equations (32.12) and (32.3), we obtain
$$\beta^{4(1)}(t,r) = -2 H^{(1)}(t,r) + \beta^{4(3)}(r).$$
Here $\beta^{4(3)}(r)$ is an arbitrary function too.
Next, equations (32.14) and (32.11) take the form:
$$e^{\beta^{4(3)}(r)}(P^{(2)})^2(K^{(1)}(r))^2\beta^{2(1)}(t,r)\beta^{3(2)}(\varphi,z)\beta^{4(2)}(\varphi,z)L_r(t,r) = -L(t,r)\beta_z^{3(2)}(\varphi,z),$$

$$e^{\beta^{3(3)}(r)}(P^{(2)})^2(K^{(1)}(r))^2\beta^{2(1)}(t,r)\beta^{3(2)}(\varphi,z)\beta^{4(2)}(\varphi,z)L_r(t,r) = -(P^{(2)})^2 L(t,r)\beta_\varphi^{4(2)}(\varphi,z).$$

Dividing left side of the first equation by the left side of the second one and the right side of the first by the right side of the second, we obtain

$$e^{-\beta^{3(3)}(r)+\beta^{4(3)}(r)} = \frac{\beta_z^{3(2)}(\varphi,z)}{(P^{(2)})^2 \beta_\varphi^{4(2)}(\varphi,z)}.$$

Clearly, both sides of this equation must be equal to the same constant value. We denote it $\ln c^{(1)}$. Then the left-hand side gives

$$\beta^{4(3)}(r) = \beta^{3(3)}(r) + c^{(1)}.$$

Equations (32.13) and (32.6) become:

$$-L(t,r)K_r^{(1)}(r) - K^{(1)}(r)L(t,r)H_r^{(1)}(t,r) + e^{2H^{(1)}(t,r)}K^{(1)}(r)(\beta^{2(1)}(t,r))^2 L_r(t,r) = 0,$$

$$L(t,r)K_r^{(1)}(r) + K^{(1)}(r)L(t,r)\beta_r^{3(3)}(r) - K^{(1)}(r)L(t,r)H_r^{(1)}(t,r) +$$
$$+ e^{2H^{(1)}(t,r)}K^{(1)}(r)(\beta^{2(1)}(t,r))^2 L_r(t,r) = 0$$

If we now subtract the second equation from the first one, we get (to an approximation of a nonzero factor)

$$2K_r^{(1)}(r) + K^{(1)}(r)\beta_r^{3(3)}(r) = 0,$$

which yields

$$\beta^{3(3)}(r) = c^{(2)} - 2\ln K^{(1)}(r),$$

where $c^{(2)} = const$.

Then the equation (20), again to an approximation of a nonzero coefficient, becomes

$$-e^{2H^{(1)}(t,r)}(K^{(1)}(r))^2 + e^{4H^{(1)}(t,r)}(K^{(1)}(r))^2 \beta^{2(1)}(t,r) + e^{2c^{(2)}}(\beta^{3(2)}(\varphi,z))^2 +$$
$$+ e^{2c^{(1)}+2c^{(2)}}(P^{(2)})^2(\beta^{4(2)}(\varphi,z))^2 = 0$$

The first two terms in the left-hand side depends on $t$ and $r$ only, and the third and fourth on $\varphi$ and $z$ only. So this equation is divided into two:

$$-e^{2H^{(1)}(t,r)}(K^{(1)}(r))^2 + e^{4H^{(1)}(t,r)}(K^{(1)}(r))^2 \beta^{2(1)}(t,r) + \lambda^2 = 0$$

and

$$e^{2c^{(2)}}(\beta^{3(2)}(\varphi,z))^2 + e^{2c^{(1)}+2c^{(2)}}(P^{(2)})^2(\beta^{4(2)}(\varphi,z))^2 - \lambda^2 = 0,$$

where $\lambda = const$.

And we get

$$\beta^{2(1)}(t,r) = e^{-2H^{(1)}(t,r)}\sqrt{e^{2H^{(1)}(t,r)} - \frac{\lambda^2}{(K^{(1)}(r))^2}}$$

and

$$\beta^{4(2)}(\varphi,z) = \frac{e^{-c^{(1)}-c^{(2)}}\sqrt{\lambda^2 - e^{2c^{(2)}}(\beta^{3(2)}(\varphi,z))^2}}{P^{(2)}}.$$

Then equation (32.14), after some transformations takes the form

$$\frac{e^{-2H^{(1)}(t,r)}P^{(2)}\sqrt{e^{2H^{(1)}(t,r)} - \frac{\lambda^2}{(K^{(1)}(r))^2}}L_r(t,r)}{L(t,r)} + \frac{\beta_z^{3(2)}(\varphi,z)}{\beta^{3(2)}(\varphi,z)\sqrt{\lambda^2 - e^{2c^{(2)}}(\beta^{3(2)}(\varphi,z))^2}} = 0.$$

For obvious reasons, this equation is divided into two. We solve the second:

$$\frac{\beta_z^{3(2)}(\varphi,z)}{\beta^{3(2)}(\varphi,z)\sqrt{\lambda^2 - e^{2c^{(2)}}(\beta^{3(2)}(\varphi,z))^2}} - \mu = 0,$$

where $\mu = const$.

It's solution is:

$$\beta^{3(2)}(\varphi,z) = \frac{4e^{\lambda(z\mu+\beta^{3(4)}(\varphi))}\lambda}{1+4e^{2(c^{(2)}+z\lambda\mu+\lambda\beta^{3(4)}(\varphi))}},$$

where $\beta^{3(4)}(\varphi)$ is an arbitrary function.

Next, equation (32.14), after some transformations, and to an approximation of a nonzero factor in the left-hand side becomes

$$e^{2H^{(1)}(t,r)}\mu L(t,r) + P^{(2)}\sqrt{e^{2H^{(1)}(t,r)} - \frac{\lambda^2}{(K^{(1)}(r))^2}} L_r(t,r) = 0. \tag{32.14a}$$

And equation (32.11) again after some transformations and to an approximation of a nonzero factor in the left-hand side becomes

$$\frac{e^{c^{(2)}+z\lambda\mu+\lambda\beta^{3(4)}(\varphi)}\beta^{3(4)}_\varphi(\varphi)}{-1+4e^{2(c^{(2)}+z\lambda\mu+\lambda\beta^{3(4)}(\varphi))}} + \frac{e^{-2H^{(1)}(t,r)}\sqrt{e^{2H^{(1)}(t,r)} - \frac{\lambda^2}{(K^{(1)}(r))^2}} L_r(t,r)}{4L(t,r)} = 0, \tag{32.11a}$$

which also splits into two. The second of them, after some transformations looks as follows:

$$-4e^{2H^{(1)}(t,r)}P^{(2)}\nu L(t,r) + P^{(2)}\sqrt{e^{2H^{(1)}(t,r)} - \frac{\lambda^2}{(K^{(1)}(r))^2}} L_r(t,r) = 0,$$

where $\nu = const$.

Comparing this equation with (32.14a), we find:

$$\nu = -\frac{\mu}{4P^{(2)}}.$$

Then the first of the equations into which the (32.11a) is splitted will be:

$$-\frac{\mu}{4P^{(2)}} + \frac{e^{c^{(2)}+z\lambda\mu+\lambda\beta^{3(4)}(\varphi)}\beta^{3(4)}_\varphi(\varphi)}{-1+4e^{2(c^{(2)}+z\lambda\mu+\lambda\beta^{3(4)}(\varphi))}} = 0.$$

It is clear that no one function $\beta^{3(4)}(\varphi)$ can make the second term on the left-hand side of this equation to be a constant value in the case when the coefficient before $z$ is not equal to zero. Therefore, we set:

$$\lambda = 0.$$

You can, of course, let $\mu = 0$, but because we are looking for at least some solution of the system, we shall not investigate this case.

Integrating this equation with $\lambda = 0$, we find:

$$\beta^{3(4)}(\varphi) = \frac{e^{-c^{(2)}}(-1+4e^{2c^{(2)}})\mu\varphi}{4P^{(2)}} + \zeta,$$

where $\zeta = const$.

Now the equation (32.4) to an approximation of a nonzero factor in the left-hand side takes the form

$$L_r(t,r) + e^{H^{(1)}(t,r)}H^{(1)}_t(t,r) = 0.$$

And the equation (32.10) to an approximation of a nonzero factor in the left-hand side takes the form

$$L(t,r)K^{(1)}_r(r) + L(t,r)K^{(1)}(r)H^{(1)}_r(t,r) + e^{H^{(1)}(t,r)}K^{(1)}(r)H^{(1)}_t(t,r) = 0.$$

If we now multiply the penultimate equation by $K^{(1)}(r)$, the last – by -1, then add these equations and integrate the final, we get:

$$L(t,r) = e^{H^{(1)}(t,r)}K^{(1)}(r)f(t),$$

where $f(t)$ is an arbitrary function.

Now the equation (32.15) to an approximation of a nonzero factor in the left-hand side becomes

$$f(t)K^{(1)}_r(r) + f(t)K^{(1)}(r)H^{(1)}_r(t,r) + H^{(1)}_t(t,r) = 0.$$

This is a non-homogeneous partial differential equation of the first order for $H^{(1)}(t,r)$. We solve it by using the characteristic system. General resolution is:

$$\Phi(H^{(1)}(t,r)+\ln K^{(1)}(r), \int f(t)dt - \int \frac{dr}{K^{(1)}(r)}) = 0,$$

where $\Phi$ is an arbitrary function of two arguments.
We single out any particular solution, for example, the following:

$$H^{(1)}(t,r) = -\ln K^{(1)}(r) + F(-\int f(t)dt + \int \frac{dr}{K^{(1)}(r)}),$$

where $F$ is an arbitrary function.
After that, the system is fully integrated.
In order to simplify the form of the solution we can make the following substitutions:
$f(t) = \tau_t(t)$,

$$K^{(1)}(r) = \frac{1}{\rho_r(r)},$$

$$F(-\int f(t)dt + \int \frac{dr}{K^{(1)}(r)}) = \ln W(-\int f(t)dt + \int \frac{dr}{K^{(1)}(r)})$$

$P^{(2)} = p = const$.
Then the solutions of the system in a cylindrically symmetric case are:

$$\vec{\beta} = (\frac{1}{W(\rho(r)-\tau(t))\tau_t(t)}, \frac{1}{W(\rho(r)-\tau(t))\rho_r(r)}, 0, 0),$$

$(g_{ij}) = diag(W^2(\rho(r)-\tau(t))\tau_t^2(t), -W^2(\rho(r)-\tau(t))\rho_r^2(r), -W^2(\rho(r)-\tau(t)), -W^2(\rho(r)-\tau(t))p^2)$.

In the case when $\rho(r) = r$ and $\tau(t) = t$ the solutions get the most simply form:

$$\vec{\beta} = (\frac{1}{W(r-t)}, \frac{1}{W(r-t)}, 0, 0),$$

$(g_{ij}) = diag(W^2(r-t), -W^2(r-t), -W^2(r-t), -W^2(r-t)p^2)$.

Riemann-Christoffel tensor of this metric is not zero.

If solve (19a), (20), (21) (i.e., in the case when the vector function is differentiable from the right and/or from the left), in the case of Minkowski space again we get a cone of light $(\beta^1)^2 + (\beta^2)^2 + (\beta^3)^2 = 1$.
Spherically symmetric solution also does not exist.
A cylindrically symmetric solution has the form

$$\vec{\beta} = \left( \frac{\sqrt{-1+W^2(f(r)-\phi(t))}}{\phi'(t)}, \frac{\sqrt{-1+W^2(f(r)-\phi(t))}}{f'(r)}, \right.$$

$$\left. \frac{\sqrt{-1+W^2(f(r)-\phi(t))}}{\sqrt{1+c^2}}, \frac{c\sqrt{-1+W^2(f(r)-\phi(t))}}{\sqrt{1+c^2}} \right)$$

$$(g_{ij}) = diag\left( \frac{\phi'^2(t)}{-1+W^2(f(r)-\phi(t))}, -\frac{W^2(f(r)-\phi(t))f'^2(r)}{-1+W^2(f(r)-\phi(t))}, -1, -1 \right).$$

The Riemann tensor of this metric is zero.
We shall not present the process of obtaining of these results.

Spherically symmetric solution of the system (13a) (i.e. when $\vec{u} \neq \vec{\beta}$) exists and is presented in [7]. The Riemann tensor of this metric is again different from zero.

## Conclusion.

As it is known, the quaternion algebra and analysis are not widely used in physics, perhaps because of the fact that they are built in a Euclidean space. In this paper we attempt to construct an alternative theory in 4-dimensional pseudo-Riemannian space, what is our space-time is (at least on a large enough scale), was undertaken.

The different types of Cauchy-Riemann equations under different assumptions about the properties of the vector field are obtained. These equations are the system of partial differential equations of the 1st order as for the components of the metric tensor so for the components of the vector field depending on the coordinates.

And if the properties of the metric tensor as physically real object are postulated, then the idea what the vector field is remains unclear. The resulting solutions at present does not clear the picture.

Also at present the formalism circumscribed here is not yet visible close analogy to the generally accepted principle of extreme action.

Will this formalism have any employment in physics or not, the further research will show.

## References.